\documentclass[11pt,a4paper]{article}
\usepackage{amssymb, amstext, amscd, amsmath, color}
\usepackage{graphicx, url, hhline}
\usepackage[left=2.50cm, right=2.00cm, top=2.30cm, bottom=2.00cm]{geometry}

\usepackage{float}

\usepackage{longtable}
\usepackage{multicol}

\newtheorem{defin}{}
\newtheorem{saetze}[defin]{}
\newtheorem{conjec}[defin]{}
\newtheorem{lemmas}[defin]{}
\newtheorem{folger}[defin]{}
\newtheorem{bemerk}[defin]{}
\newtheorem{prop}{Proposition} 

\newcommand{\fillbox}{\mbox{$\bullet$}}

\newcommand{\Z}{\mathbb Z}

\newcommand{\R}{\mathbb R}
\newcommand{\GAP}{\mathsf{GAP}}

\begin{document}

\bigskip
\bigskip

\centerline{\large \textbf{On colourings of cubic lattices}}

\bigskip
\bigskip

\centerline{\large Igor A. Baburin}
\centerline{\emph{Email: baburinssu@gmail.com}}

\bigskip

\noindent Given the integral lattice $\Lambda^d$ in $d$-dimensional Euclidean space, partitions of the lattice nodes into orbits of finite-index subgroups of $Aut(\Lambda^d)$ have been computed for $d \leq 4$. These partitions can be interpreted as colourings of orbits defined up to permutation of colours. Complete results are obtained for $d=2$ up to 64 orbits, for $d=3$ up to 8 orbits, and for 2 orbits in dimension~4. The automorphism groups of the partitions are also determined. Our results for two orbits in dimension~3 correct the old result of H. Heesch [\emph{Z. Kristallogr.}, (1933), 85, 335--344] who overlooked one partition. 

\section{Introduction}

In 1933 \emph{Kristallmathematiker} H. Heesch published a paper \emph{Zur Topologie parallelepipedischer Gitter} where he -- by using a combinatorial approach -- constructed (in modern terms) partitions of the nodes of the square lattice ($\Lambda^2$) and the three-dimensional cubic lattice ($\Lambda^3$) into two orbits \cite{Hee33}. The orbits refer to the action of finite-index subgroups of the automorphism group of the respective lattice.\footnote{In contrast to the usual practice in number theory, here the automorphism group of a lattice is understood in a broader sense: it is generated by certain elements of $\mathrm{GL}(d, \Z)$ \emph{and} integral translations.} For illustrations he used tessellations of $\R^2$ by squares resp. $\R^3$ by cubes, and therefore he invented the natural term \emph{Zweiteilungen} for these partitions.  For the square lattice he obtained \emph{nine} partitions shown in Figure~\ref{f:sql2} below. This result was checked a number of times (see \emph{e.g.} \cite{Puz08}). In dimension 3 Heesch found 24 partitions, and this result seems to have never been checked or reproduced. In this paper, using the methods of computational group theory, we show that there are precisely 25 two-orbit partitions of the nodes of the three-dimensional cubic lattice, and go much further with the number of orbits as well as in the dimension.

\begin{center}
\begin{figure}[h]
\centering
\includegraphics{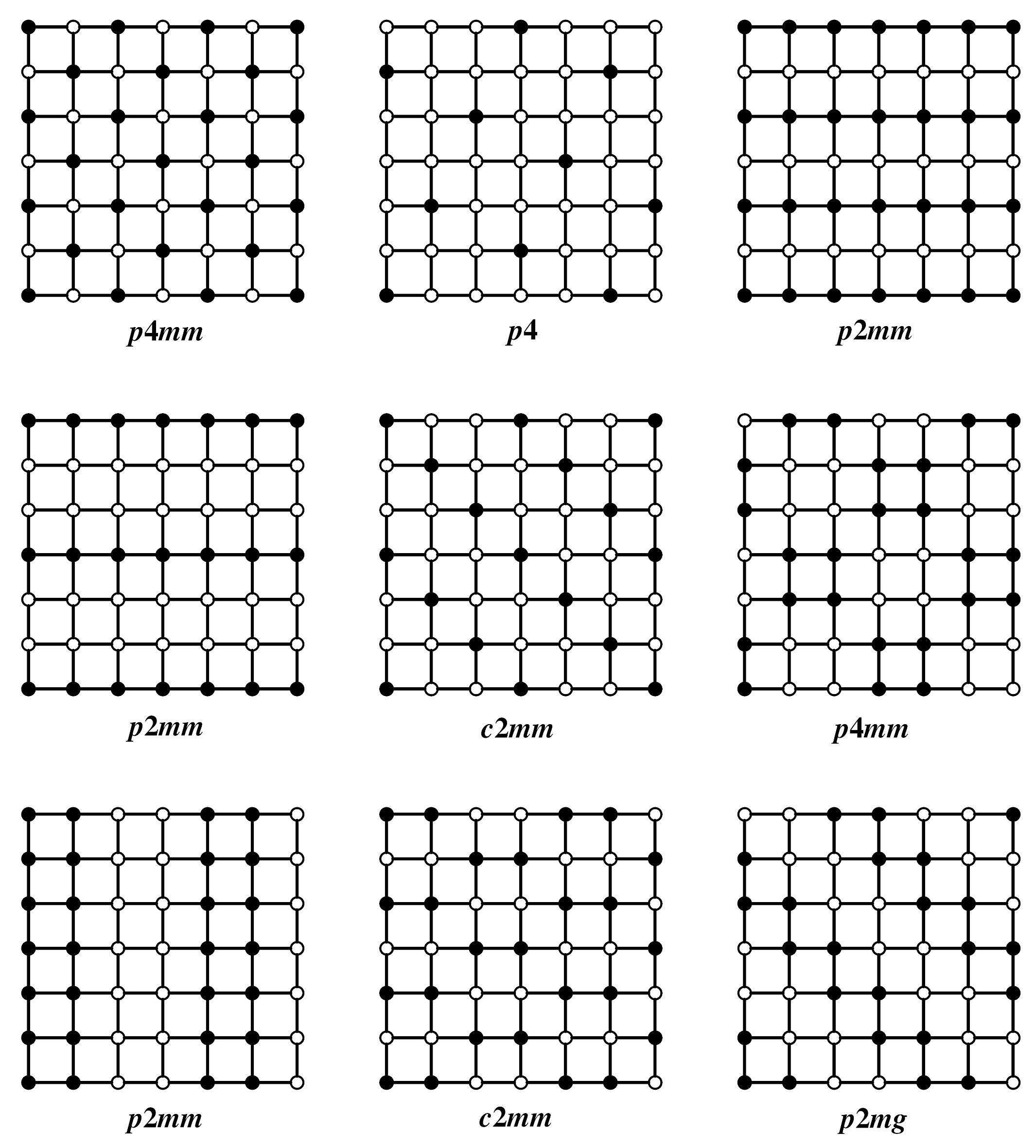}
\caption{Two-orbit partitions of the square lattice.}
\label{f:sql2}
\vspace{-2em}
\end{figure}
\end{center}

\section{Theory and computational methods}

Consider a finite set $X$, on which a group $G$ acts transitively (so $X$ is an orbit of $G$). The stabilizer of $x \in X$ in $G$ is defined as usual as $Stab_G(x) = \{ g | xg = x, g \in G \}$. We need the following elementary proposition.

\begin{prop}
\label{prop1}
Let $X$ be an orbit of a finite group $G$, and let $H < G$ with index $i = |G:H|$. \\ 
By the orbit-stabilizer lemma we have $|X| = |G|/|Stab_{G}(x)|$. Under the action of $H$, $X$ splits into $n$ orbits $X_1, ... , X_n$. Then $|X| = \sum_{k=1}^{n} |H|/|Stab_{H}(x_k)|$. By equating the two expressions for $|X|$ we obtain $i = \sum_{k=1}^{n} |Stab_{G}(x) : Stab_{H}(x_k)|$.
\end{prop}

\noindent This proposition carries over to crystallographic groups $G$ and $H$, both acting on a discrete point set $X \subset \R^d$, in which case orbit cardinalities are counted modulo translation lattice of $H$ \cite{Won93}. We remind here that a crystallographic group $G$ of $\R^d$ is an extension of $T \cong \Z^d$ by a finite group $P = G/T$ so that $P$ acts faithfully on the lattice $T$. The group $T$ is the translation subgroup of $G$ that we shall denote by $T(G)$, and $P$ is its point group.

Let $X$ be an orbit of a crystallographic group in $\R^d$, and let $Aut(X)$ be the group of all its isometries. Partitions of $X$ into $n$ orbits can be determined in the following steps:

\begin{enumerate}
\setlength\itemsep{-0.02in}
\item[(a)] Compute subgroups of $Aut(X)$ with index at most $n \times |Stab_{Aut(X)}(x)|$ (up to conjugacy in $Aut(X)$). Candidates for subgroup indices can be found by running over all possible $n$-tuples of divisors of $|Stab_{Aut(X)}(x)|$.
\item[(b)] Compute orbits of subgroups on $X$.
\item[(c)] Order by inclusion subgroups with the same number of orbits on $X$ (replacing a given subgroup by a conjugate, if necessary).
\end{enumerate}

\noindent This computational scheme has been implemented in the computer algebra system $\GAP$ \cite{GAP} using different methods that are described below.

Step (a) is best carried out using a convenient finite presentation for $Aut(X)$. For cubic lattices ($X = \Lambda^d$) these are Coxeter-type presentations (Table~\ref{t:table1}) from which it is not so straightforward to see that $Aut(\Lambda^d) \cong \Z^d \rtimes ( S_d \rtimes (C_2)^d )$. Since the groups in Table~\ref{t:table1} are polycyclic, the algorithms for polycyclic groups can be employed as well (to do so, polycyclic presentations are necessary which are different from those in Table~\ref{t:table1}). 

Let us fix subgroup $H \leq Aut(X)$ of finite index. $H$ is either transitive on $X$ or defines a partition $\Pi$ of $X$ into two or more orbits. To compute orbits of $H$ on $X$ at Step~(b), the Schreier--Sims method can be applied either for finite matrix groups (orbits are determined mod~$T(H)$) or finite permutation groups. The finite permutation group arises as the image of the action of $H$ on~$X$ mod $T(H)$ \cite{Felix11}. Our own implementation of the Schreier--Sims method (following~\cite{Cannon84}) for matrix groups works slightly faster in dimensions $d \leq 4$ than the built-in $\GAP$ orbit functions for permutation groups.  

Step (c) that is performed with the usual tools of computational group theory directly yields the full automorphism groups of the partitions. The automorphism group of the partition, $Aut(\Pi)$, is understood here as the group \emph{fixing the orbits} but \emph{not} permuting them (it coincides with the so-called \emph{colour fixing group} from \cite{Felix11}). The automorphisms so defined are contained in the intersection of all orbit stabilizers in $Aut(X)$ -- the point of view originally followed by H.~Heesch \cite{Hee33}.

In dimensions 2 and 3 the presentations from Table~\ref{t:table1} allow to compute subgroups up to very high indices. Using the $\GAP$ tools for finitely-presented groups, it is possible to go, for example, up to index 512 for $d=2$, and up to index 384 for $d=3$ with a moderate computational effort. These subgroup data are sufficient for partitioning the nodes of the square lattice in up to 64 orbits ($512 = 8 \times 64$), and for the three-dimensional cubic lattice up to 8 orbits ($384 = 48 \times 8$).

For $d=4$ this brute-force approach is not appropriate anymore because the number of subgroups to be considered already for two-orbit partitions becomes overwhelming (in the order of a million or so). In this case we adopted a two-step procedure:


\begin{enumerate}
\setlength\itemsep{-0.02in}
\item[(a)] compute node-transitive subgroups of $Aut(\Lambda^4)$, and
\item[(b)] for each node-transitive subgroup compute its \emph{maximal subgroups},
\end{enumerate}

\noindent from which we filtered out those acting on $\Lambda^4$ with two orbits. 

The indices of node-transitive subgroups of $Aut(\Lambda^4)$ correspond to divisors of $384(=2^4 \times 4!)$. When applying algorithms for finitely-presented groups, the expressions for generators of subgroups of high indices can be become rather long and practically unmanageable. Therefore, we have chosen a combined strategy: subgroups of $Aut(\Lambda^4)$ with indices up to 64 were computed from the presentation in Table~\ref{t:table1}. Subgroups with indices $\{96, 128, 192, 384\}$ were determined using algorithms for polycyclic groups \cite{Polyc2020}. Subgroups of index 128 were computed with both methods for checking purposes.

{\tiny
\begin{table}[!h]
\begin{center}
\caption{Presentations of automorphism groups for cubic lattices in $\R^d$ ($d=2, 3, 4$).}
\label{t:table1}
\begin{tabular}{|l|l|l|}
\hline
$d$ & Generators & Presentation \\ \hline
2 & $a=m_{0y}: -x, y$ & $\langle a, b, c | a^2, b^2, c^2, (ba)^2, (cb)^4, (ca)^4 \rangle$ \\
  & $b=m_{x\frac{1}{2}}: x, 1-y$ & \\
  & $c=m_{xx}: y, x$ &     \\
\hline
3 & $a=m_{xy0}: x, y, -z$ & $\langle a, b, c, d | a^2, b^2, c^2, d^2, (ac)^2, (bd)^2, (ba)^2, (cd)^3, (da)^4, (cb)^4 \rangle$ \\
  & $b=m_{\frac{1}{2}yz}: 1-x, y, z$ & \\
  & $c=m_{xxz}: y, x, z$  & \\
  & $d=m_{xyy}: x, z, y$ & \\
\hline
4 & $a=m_{xxzw}: y, x, z, w$ & $\langle a, b, c, d, f | a^2, b^2, c^2, d^2, f^2, (ba)^2, (cd)^2, (cb)^2, (af)^2, (cf)^2$, \\
  & $b=m_{xy0w}: x, y, -z, w$ & $  \hspace{2cm} (df)^3, (da)^3, (dfbf)^2, (ca)^4, (db)^4, (bf)^4 \rangle$ \\
  & $c=m_{\frac{1}{2}yzw}: 1-x, y, z, w$ & \\
  & $d=m_{xyyw}: x, z, y, w$ & \\
  & $f=m_{xyww}: x, y, w, z$ & \\
\hline
\end{tabular}
\end{center}
\end{table}
}

The computation of maximal subgroups for node-transitive groups of $\Lambda^4$ was performed as described by Eick \emph{et al.} (1997) \cite{Eick97} and implemented in the \emph{Cryst} package of the $\GAP$ system \cite{Cryst2019}. The index of maximal subgroups can be effectively bounded from Proposition~\ref{prop1} in combination with the arguments making use of the subgroup structure of crystallographic groups. The number of \emph{translationengleiche} subgroups is always finite whereas the indices of maximal \emph{klassengleiche} subgroups are restricted to powers of $2, 3, 5, 7$ in the current context. Moreover, subgroup indices can be bounded depending on the node stabilizer: for groups with node stabilizers of order 2 it is sufficient to compute maximal subgroups up to index 4 whereas in the case of trivial stabilizers only subgroups of index~2 are relevant. Note that this procedure can be generalised in a straightforward way for computing partitions into three or more orbits (that naturally implies higher bounds for the indices of maximal subgroups). For example for three orbits, maximal subgroups of both node-transitive groups of~$\Lambda^4$ and \emph{all} the groups (not only the full automorphism groups) of the two-orbit partitions have to be considered. However, with increasing number of orbits the method quickly becomes inefficient due to a high number of conjugacy checks (w.r.t. $Aut(\Lambda^4)$) to be performed.

Let us illustrate this procedure by working out two-orbit partitions for the linear lattice $\Lambda^1$. We have: $Aut(\Lambda^1) = \langle a, b \rangle$, where $a$ stands for the translation, $b$ -- for the inversion through the node at~$\bf 0$. Node-transitive groups of $\Lambda^1$ are $\langle a, b \rangle, \langle a \rangle, \langle a^2, ba \rangle$. From Proposition~\ref{prop1} it follows that two-orbit partitions arise for subgroups of $Aut(\Lambda^1)$ with index at most 4 which can in turn be obtained from maximal subgroups of node-transitive groups with index at most~3 (Table~\ref{t:table2}).

{\tiny
\begin{table}[ht]
\vspace{-1em}
\begin{center}
\caption{Two-orbit partitions of $\Lambda^1$ from its node-transitive groups.}
\label{t:table2}
\begin{tabular}{|l|l|l|}
\hline
Group & Subgroups & Orbit partition \\ \hline
$\langle a, b \rangle$ & $\langle a \rangle, \langle a^2, ba \rangle$ & node-transitive\\
& $\langle a^2, b \rangle$ & $...+-...$ \\
& $\langle a^3, b \rangle$ & $...++-...$ \\
\hline
$\langle a \rangle$ &  $\langle a^2 \rangle$ & $...+-...$ \\
\hline
$\langle a^2, ba \rangle$ & $\langle a^4, ba \rangle$ & $...++--...$ \\
& $\langle a^2 \rangle$ &  $...+-...$ \\
\hline
\end{tabular}
\end{center}
\vspace{-1em}
\end{table}
}

Let us emphasize that the employed group-theoretic tools are applicable to compute orbit partitions for any periodic graph (\emph{cf.} \cite{Jun23}). Owing to high symmetry, cubic integral lattices stand out as computationally most demanding case, and therefore represent a good benchmark for computational methods.

\section{Results} 

Table~\ref{t:table3} serves as an overview of the obtained results. Therein, for each lattice $\Lambda^d$ in dimension $d=2,3,4$ a sequence of numbers is given whose $k$th term corresponds to the number of partitions into $(k+1)$ orbits. Keeping in mind applications to structural chemistry, for $\Lambda^3$ we provide (in round brackets) the number of partitions where only nodes from different orbits are adjacent. As already noticed by Heesch \cite{Hee33}, among orbit partitions of $\Lambda^d$ there always exist those which are derived from partitions of $\Lambda^{d-1}$ by a direct superposition in the $d$th dimension. For example, Figure~\ref{f:sql2} clearly shows that partitions with symmetry $p2mm$ are derived from the two-orbit partitions of the linear lattice $\Lambda^1$. This property (among other things) has been useful for checking purposes.

{\tiny
\begin{table}[!h]
\vspace{-1em}
\begin{center}
\caption{Orbit partitions of cubic lattices.}
\label{t:table3}
\begin{tabular}{|l|l|l|}
\hline
$d$ & \parbox[t][1cm]{3.8cm}{\raggedright No. of node-transitive groups of $\Lambda^d$} & No. of partitions \\ \hline
2 & 36 & 9, 22, 44, 39, 80, 47, 96, 81, 104, 65, 157, 75, 119, 129, \\
  & & 160, 78, 196, 83, 201, 153, 147, 111, 255, 148, 170, 160, \\
  & & 244, 116, 287, 125, 257, 188, 191, 192, 362, 138, 215, 207, \\
  & & 331, 160, 334, 176, 281, 287, 258, 157, 422, 203, 341, 235, \\
  & & 316, 208, 367, 261, 413, 250, 288, 201, 542, 223, 278, 348, 409  \\
\hline
        3 & 786 & 25(1), 80(4), 275(34), 281(17), 800(103), 432(54), 1469(282) \\
\hline
4 & 38725 & 73 \\
\hline
\end{tabular}
\end{center}
\end{table}
}

\begin{center}
\begin{figure}[!h]
\vspace{-1em}
\centering
\includegraphics{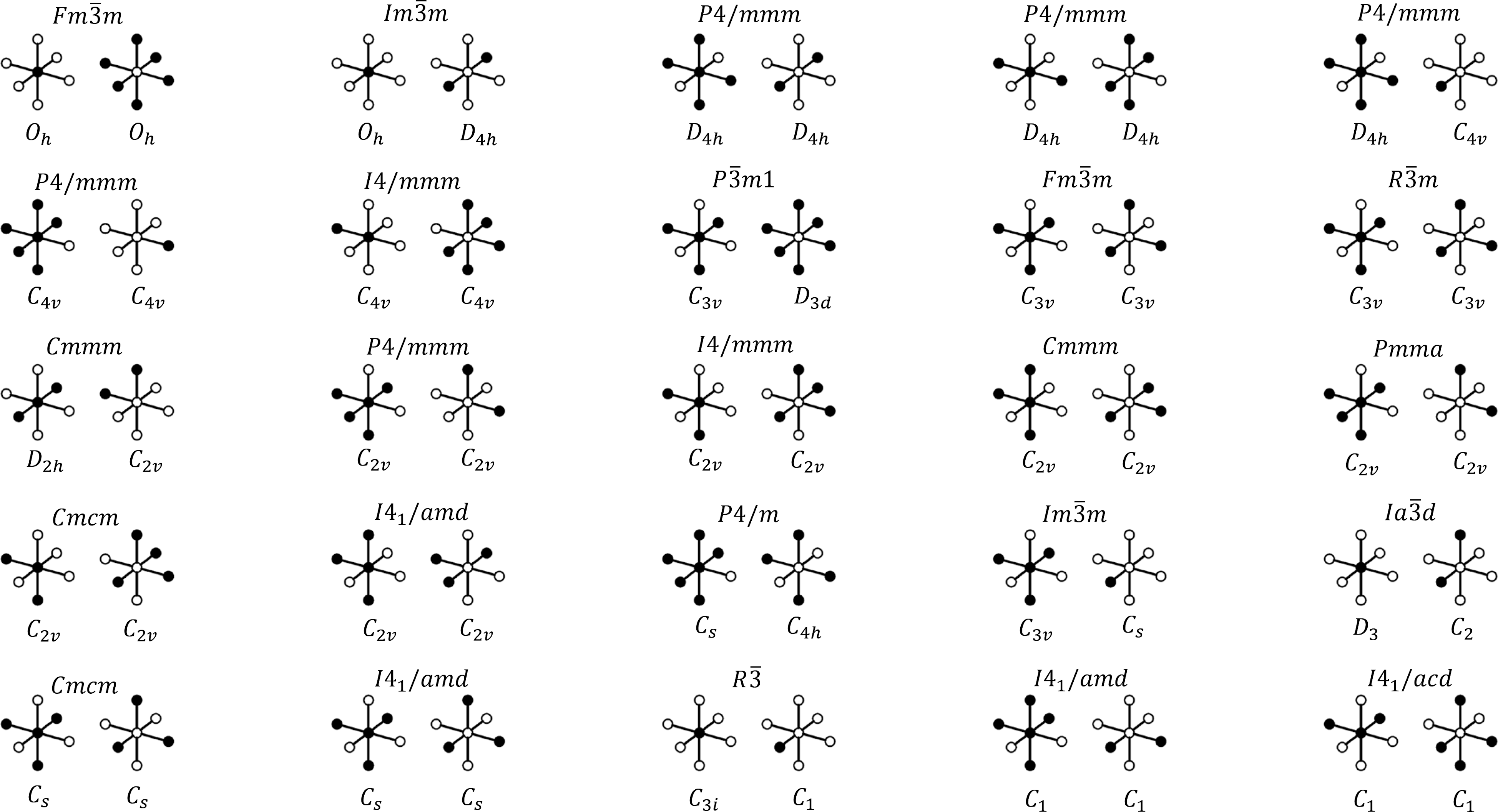}
\caption{Node neighbourhoods in the two-orbit partitions of $\Lambda^3$.}
\label{f:pcu2}
\vspace{-2em}
\end{figure}
\end{center}

{\tiny
\begin{table}[!ht]
\begin{center}
\caption{Two-orbit partitions of $\Lambda^3$.}
\label{t:table4}
\begin{tabular}{|l|l|l|l|l|l|l|l|}
\hline
No. & $Aut(\Pi)$ & \parbox[t]{1cm}{Index \\$i_t \cdot i_k$} & \parbox[t]{2.3cm}{Wyckoff pos. ${\bf{a}}^2, {\bf{b}}^2, {\bf{c}}^2$} & No. & $Aut(\Pi)$ & \parbox[t]{1cm}{Index \\$i_t \cdot i_k$} & \parbox[t]{2.3cm}{Wyckoff pos. ${\bf{a}}^2, {\bf{b}}^2, {\bf{c}}^2$} \\
\hline
1 & $Fm\bar3m$ & $2 \cdot 1$ & $4a~0, 0, 0$ & 13 & $I4/mmm$ & $8 \cdot 3$ & $8i~0, 3/4, 0$  \\
   &    &    &   $4b~0, 0, 1/2$ &  &    &    &   $8j~0, 3/4, 1/2$ \\
   &    &    &   16, 16, 16 &  &    &    &   32, 32, 16  \\
\hline
2 & $Im\bar3m$ & $4 \cdot 1$ & $2a~0, 0, 0$ & $14^\ast$ & $Cmmm$ & $4 \cdot 6$ & $4h~3/8, 1/2, 1/2$  \\
   &    &    &   $6b~0, 0, 1/2$ &  &    &    &   $4h~3/8, 0, 1/2$ \\
   &    &    &   16, 16, 16 &  &    &    &   64, 16, 4  \\
\hline
$3^\ast$ & $P4/mmm$ & $2 \cdot 3$ & $1c~1/2, 1/2, 0$ & $15^\ast$ & $Pmma$ & $4 \cdot 6$ & $2e~1/4, 0, 7/8$  \\
   &    &    &   $1d~1/2, 1/2, 1/2$ &  &    &    &   $2e~1/4, 0, 3/8$ \\
   &    &    &   4, 4, 16 &  &    &    &   8, 4, 32  \\
\hline
$4^\ast$ & $P4/mmm$ & $2 \cdot 3$ & $1b~0, 0, 1/2$ & 16 & $Cmcm$ & $4 \cdot 6$ & $4c~1/2, 3/8, 3/4$  \\
   &    &    &   $1d~1/2, 1/2, 1/2$ &  &    &    &   $4c~0, 3/8, 3/4$ \\
   &    &    &   8, 8, 4 &  &    &    &   16, 32, 8  \\
\hline
$5^\ast$ & $P4/mmm$ & $3 \cdot 3$ & $1c~1/2, 1/2, 0$ & 17 & $I4_1/amd$ & $8 \cdot 3$ & $8e~1/2, 1/4, 0$  \\
   &    &    &   $2h~1/2, 1/2, 1/3$ &  &    &    &   $8e~0, 1/4, 0$ \\
   &    &    &   4, 4, 36 &  &    &    &   16, 16, 64  \\
\hline
$6^\ast$ & $P4/mmm$ & $4 \cdot 3$ & $2g~0, 0, 7/8$ & $18^\ast$ & $P4/m$ & $5 \cdot 6$ & $4j~1/5, 3/5, 0$  \\
   &    &    &   $2g~0, 0, 3/8$ &  &    &    &   $1a~0, 0, 0$ \\
   &    &    &   4, 4, 64 &  &    &    &   20, 20, 4  \\
\hline
7 & $I4/mmm$ & $4 \cdot 3$ & $4e~1/2, 1/2, 3/8$ & 19 & $Im\bar3m$ & $32 \cdot 1$ & $16f~3/8, 3/8, 3/8$  \\
   &    &    &   $4e~0, 0, 3/8$ &  &    &    &   $48k~3/8, 3/8, 7/8$ \\
   &    &    &   8, 8, 64 &  &    &    &   64, 64, 64  \\
\hline
8 & $P\bar3m1$ & $3 \cdot 4$ & $2d~1/3, 2/3, 1/3$ & 20 & $Ia\bar3d$ & $32 \cdot 1$ & $16b~3/8, 5/8, 1/8$  \\
   &    &    &   $1a~0, 0, 0$ &  &    &    &   $48g~3/8, 5/8, 3/8$ \\
   &    &    &   8, 8, 12 &  &    &    &   64, 64, 64  \\
\hline
9 & $Fm\bar3m$ & $16 \cdot 1$ & $32f~3/8, 3/8, 7/8$ & 21 & $Cmcm$ & $8 \cdot 6$ & $8g~3/8, 1/8, 3/4$  \\
   &    &    &   $32f~3/8, 3/8, 3/8$ &  &    &    &   $8g~3/8, 5/8, 3/4$ \\
   &    &    &   64, 64, 64 &  &    &    &   64, 32, 8  \\
\hline
10 & $R\bar3m$ & $4 \cdot 4$ & $6c~0, 0, 3/8$ & 22 & $I4_1/amd$ & $16 \cdot 3$ & $16h~0, 0, 3/8$  \\
   &    &    &   $6c~0, 0, 7/8$ &  &    &    &   $16h~0, 0, 7/8$ \\
   &    &    &   8, 8, 192 &  &    &    &   32, 32, 64  \\
\hline
$11^\ast$ & $Cmmm$ & $3 \cdot 6$ & $2b~0, 1/2, 0$ & 23 & $R\bar3$ & $7 \cdot 8$ & $3a~0, 0, 0$  \\
   &    &    &   $4g~1/6, 0, 0$ &  &    &    &   $18f~5/21, 4/21, 1/3$ \\
   &    &    &   72, 8, 4 &  &    &    &   56, 56, 12  \\
\hline
$12^\ast$ & $P4/mmm$ & $8 \cdot 3$ & $4l~3/4, 0, 0$ & 24 & $I4_1/amd$ & $32 \cdot 3$ & $32i~3/8, 1/8, 0$  \\
   &    &    &   $4n~1/4, 1/2, 0$ &  &    &    &   $32i~3/8, 1/8, 1/2$ \\
   &    &    &   32, 32, 4 &  &    &    &   64, 64, 64  \\
\hline
25 & $I4_1/acd$ & $32 \cdot 3$ & $32g~1/8, 3/8, 0$ & \multicolumn{3}{l}{} \\
   &    &    &  $32g~1/8, 3/8, 1/2$ & \multicolumn{3}{l}{} \\
   &    &    &  64, 64, 64  & \multicolumn{3}{l}{} \\
\cline{1-4}
\end{tabular}
\end{center}
\end{table}
}

\newpage

Table~\ref{t:table4} in combination with Figure~\ref{f:pcu2} characterises the two-orbit partitions of $\Lambda^3$. For each partition $\Pi$, its automorphism group, $Aut(\Pi)$, is given together with the coordinate description (Wyckoff positions and the diagonal components of the metric tensor) respecting the conventions of the International tables for crystallography \cite{IT}. Origin is always assumed at the center of symmetry. Metrical parameters refer to nearest-neighbour distance~=~2. Partitions are ordered by indices of $Aut(\Pi)$ in $Aut(\Lambda^3)$. Each index~$i$ is decomposed into two multipliers $i = i_ti_k$: the first one defines the unit-cell enlargement (the lowering of translational symmetry) whereas the second one indicates the reduction of the point-group (orthogonal) symmetry. Partitions that can be derived by a direct superposition of those of the square lattice are marked by an asterisk. Figure~\ref{f:pcu2} shows node neighbourhoods together with node stabilizers (in the crystallographic Schoenflies notation) in $Aut(\Pi)$. It can be noticed that some neighbourhoods show higher local symmetry than allowed by the stabilizers (for example, this is always the case for nodes with trivial stabilizers). Moreover, the following partitions have identical neighbourhoods:
\centerline{$\{14, 17\}, \{12, 15, 24\}, \{13, 16, 25\}, \{9, 10, 21, 22\}$.}
\noindent Considering neighbourhoods of radius 2 makes the differences apparent. 

\begin{center}
\begin{figure}[h]
\centering
\includegraphics{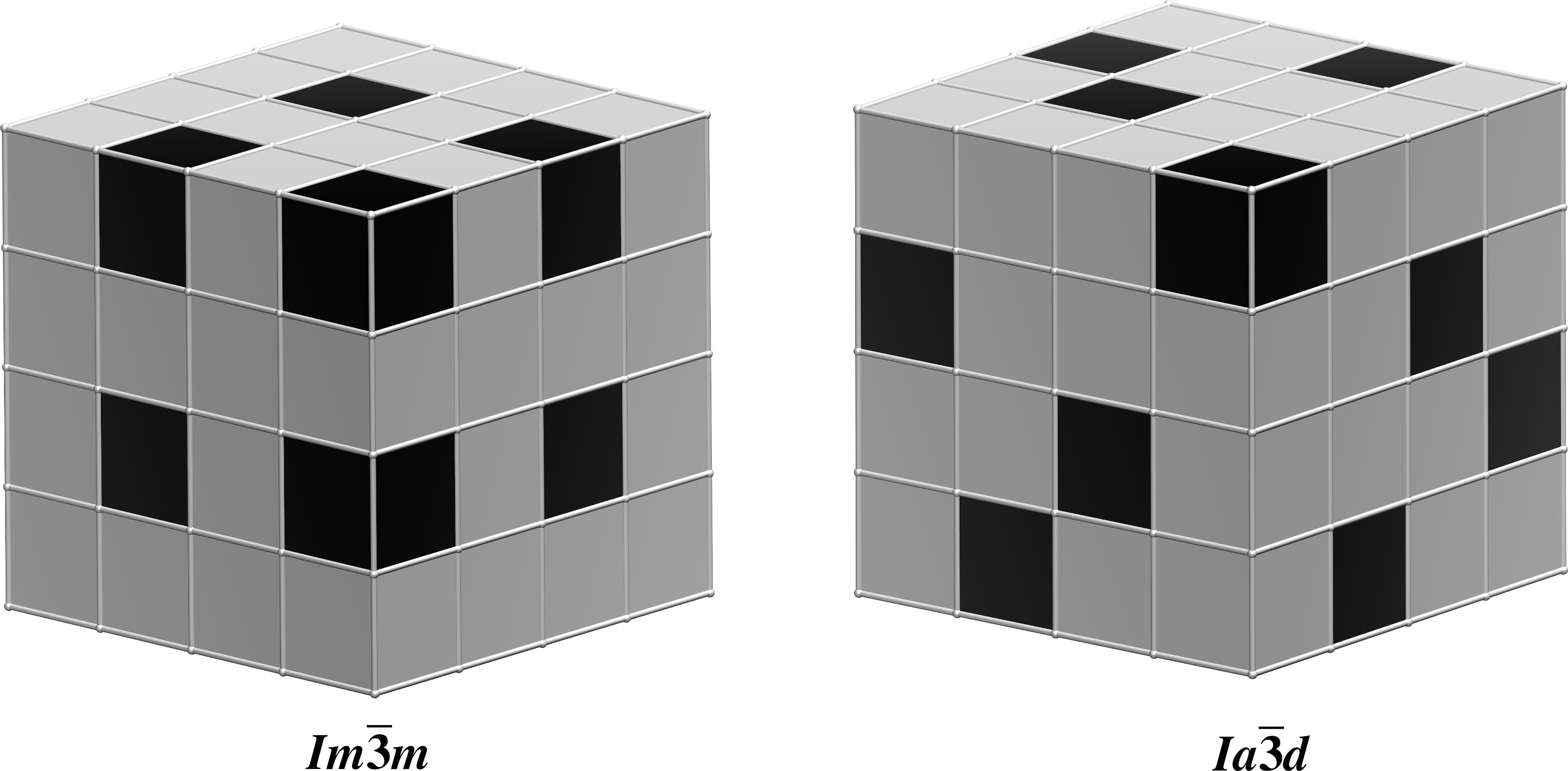}
\caption{The two-orbit partitions of $\Lambda^3$ with symmetry $Im\bar3m$ (No.~2) and $Ia\bar3d$ (No.~20). Lattice nodes are drawn as cubes.}
\label{f:Ia3d}
\vspace{-1em}
\end{figure}
\end{center}

The partition overlooked by Heesch corresponds to No.~20. Partitions No.~2 and No.~20 are compared in Figure~\ref{f:Ia3d} where three-fold mirror symmetry manifests itself in No.~2, and only three-fold rotational symmetry in No.~20. Moreover, if both partitions are referred to the same space group $Ia\bar3d$ (that implies a symmetry reduction for the partition No.~2), their duality comes to light: the vertices of cubes in No.~20 correspond to the centroids of cubes in No.~2 (it is hard to illustrate satisfactorily but can be verified by computation). This special property might explain the omission in the list of Heesch. This example nevertheless demonstrates that there is always something insightful behind the mistakes of the people from the past.

The results for $\Lambda^4$ are presented in a more condensed form due to a much larger amount of the material. Figure~\ref{f:pcu4d} shows all combinations of node neighbourhoods occurring in the two-orbit partitions of $\Lambda^4$. For each neighbourhood type, in Table~\ref{t:table5} the automorphism groups of the partitions are given following the notation of \cite{BBNWZ} for four-dimensional crystallographic groups, and the structure of node stabilizers. The 25 partitions of $\Lambda^4$ originating from a direct superposition of those of $\Lambda^3$ are marked by asterisks (their number 25 accidentally coincides with the number of neighbourhood types in Figure~\ref{f:pcu4d}), and are collected in Table~\ref{t:table6}. Various projections and applications of the two-orbit partitions for $\Lambda^4$ remain to be explored in the future.

The data on partitions of $\Lambda^2$ and $\Lambda^3$ into three and more orbits are given in the Supporting information in a form similar to that of Table~\ref{t:table4}.

{\tiny
\begin{table}[!ht]
\begin{center}
\caption{Two-orbit partitions of $\Lambda^4$.}
\label{t:table5}
\small Notation for four-dimensional space groups from \cite{BBNWZ}.
\begin{tabular}{|l|l|}
\hline
Type & $Aut(\Pi)$ \{Stabilizers\} \\ \hline
I & 32/21/2/1 $\{S_4 \rtimes (C_2)^4; S_4 \rtimes (C_2)^4\}$ \\
\hline
II & 32/11/1/1 $\{GL(2,3); S_3\}$  \\
\hline
III & 18/2/1/1 $\{(C_4 \times C_2) : C_2; C_2 \times C_2\}$  \\
\hline
IV & 24/5/1/1 $\{C_2 \times S_4; S_4\}$  \\
\hline
V & 25/11/5/1 $\{C_2 \times S_4; C_2 \times S_4\}$  \\
\hline
VI & (*) 25/11/4/1 $\{C_2 \times C_2 \times S_4; C_2 \times C_2 \times S_4\}$  \\
\hline
VII & (*) 25/11/2/1 $\{C_2 \times C_2 \times S_4; C_2 \times C_2 \times D_4\}$  \\
\hline
VIII & (*) 14/4/1/1 $\{C_6 \times C_2; C_2\}$  \\
\hline
IX & (*) 25/11/2/3 $\{D_6; C_2 \times C_2\}$  \\
\hline
X & (*) 14/10/3/1 $\{C_2 \times D_6; D_6\}$  \\
\hline
XI & (*) 13/10/3/1 $\{C_2 \times D_4; C_2 \times D_4\}$; 13/10/5/3 $\{(C_2)^3; (C_2)^3\}$ \\
 & 19/6/2/6 $\{C_2 \times C_2; C_2 \times C_2\}$  \\
\hline
XII & (*) 19/6/1/1 $\{D_4 \times D_4; D_4 \times D_4\}$; 32/17/2/3 $\{(C_2)^4 : C_2; (C_2)^4 : C_2\}$  \\
\hline
XIII & (*) 19/3/1/1 $\{C_4 \times D_4; D_4\}$  \\
\hline
XIV & (*) 13/10/2/1 $\{C_2 \times C_2 \times D_4; C_2 \times D_4\}$  \\
\hline
XV & (*) 13/10/2/1 $\{C_2 \times D_4; C_2 \times D_4\}$; (*) 13/10/3/33 $\{(C_2)^3; (C_2)^3\}$ \\
 & 19/6/2/3 $\{C_2 \times C_2; C_2 \times C_2\}$  \\
\hline
XVI & (*) 25/11/1/1 $\{C_2 \times C_2 \times S_4; C_2 \times C_2 \times S_4\}$  \\
\hline
XVII & (*) 25/11/1/1 $\{C_2 \times C_2 \times S_4; C_2 \times S_4\}$  \\
\hline
XVIII & (*) 25/11/1/1 $\{C_2 \times S_4; C_2 \times S_4\}$  \\
\hline
XIX & (*) 25/11/2/1 $\{D_6; C_2 \times C_2\}$  \\
\hline
XX & (*) 19/6/1/1 $\{C_2 \times D_4; C_2 \times D_4\}$; (*) 13/10/1/22 $\{C_2 \times D_4; C_2 \times D_4\}$ \\
 & 32/17/2/3 $\{C_2; C_2\}$; (*) 13/10/3/33 $\{C_2; C_2\}$ \\
\hline
XXI & (*) 13/10/3/1 $\{(C_2)^3; (C_2)^3\}$; 13/10/3/43 $\{(C_2)^3; (C_2)^3\}$ \\
 & 18/5/1/1 $\{(C_2)^3; (C_2)^3\}$; (*) 6/3/2/3 $\{(C_2)^3; (C_2)^3\}$ \\
 & 13/10/5/3 $\{C_2 \times C_2; C_2 \times C_2\}$; 13/10/5/6 $\{C_2 \times C_2; C_2 \times C_2\}$ \\
 & 32/13/2/5 $\{C_2; C_2\}$; (*) 13/10/3/37 $\{C_2; C_2\}$ \\
 & 13/10/5/3 $\{C_2; C_2\}$; 18/2/4/3 $\{C_2; C_2\}$ \\
 & 19/6/2/6 $\{C_1; C_1\}$; 32/13/3/5 $\{C_1; C_1\}$ \\
 & 18/2/4/3 $\{C_1; C_1\}$; 18/3/5/2 $\{C_1; C_1\}$ \\
\hline
XXII & (*) 25/11/4/1 $\{D_6; D_6\}$; (*) 14/10/1/1 $\{D_6; D_6\}$ \\
 & (*) 13/10/3/33 $\{C_2 \times C_2; C_2 \times C_2\}$; (*) 6/3/2/3 $\{C_2 \times C_2; C_2 \times C_2\}$ \\
 & 13/10/3/43 $\{C_2; C_2\}$; 18/5/1/1 $\{C_2; C_2\}$ \\
 & 13/10/5/3 $\{C_1; C_1\}$; 13/10/5/6 $\{C_1; C_1\}$ \\
\hline
XXIII & 19/6/2/1 $\{C_2 \times D_4; C_2 \times D_4\}$; 13/10/3/12 $\{C_2 \times D_4; C_2 \times D_4\}$ \\
 & 32/17/2/8 $\{C_2; C_2\}$; 13/10/5/4 $\{C_2; C_2\}$ \\
\hline
XXIV & 25/11/5/1 $\{D_6; D_6\}$; 14/10/2/1 $\{D_6; D_6\}$ \\
 & 13/10/5/3 $\{C_2 \times C_2; C_2 \times C_2\}$; 6/3/6/2 $\{C_2 \times C_2; C_2 \times C_2\}$ \\
 & 13/10/3/47 $\{C_2; C_2\}$; 18/5/1/2 $\{C_2; C_2\}$ \\
 & 13/10/5/4 $\{C_1; C_1\}$; 13/10/5/5 $\{C_1; C_1\}$ \\
\hline
XXV & 32/21/2/1 $\{S_4; S_4\}$; 24/5/5/2 $\{S_4; S_4\}$ \\
 & 32/17/2/3 $\{D_4; D_4\}$; 18/5/3/13 $\{D_4; D_4\}$ \\
 & 24/5/1/1 $\{S_3; S_3\}$; 14/10/2/1 $\{S_3; S_3\}$ \\
 & 13/10/3/12 $\{C_2 \times C_2; C_2 \times C_2\}$; 13/10/3/34 $\{C_2 \times C_2; C_2 \times C_2\}$ \\
 & 13/10/5/3 $\{C_2; C_2\}$; 18/5/5/2 $\{C_2; C_2\}$ \\
 & 6/3/6/11 $\{C_2; C_2\}$  \\
\hline
\end{tabular}
\end{center}
\end{table}
}

\section{Comments on earlier works}

Theory of colour symmetry received much attention and development in the 1950--1980s when it was established as an interdisciplinary direction of mathematical crystallography with applications mainly to the description of magnetic structures (for an overview see \emph{e.g.} \cite{Sen88}). However, colourings considered at that time were of different kind compared to the subject of this paper: given a crystal pattern with symmetry $G$ and its subgroup $H$ of finite index $n$, the pattern was coloured by $n$ colours corresponding to the cosets of $H$ in $G$. 


\begin{center}
\begin{figure}[!ht]
\centering
\includegraphics{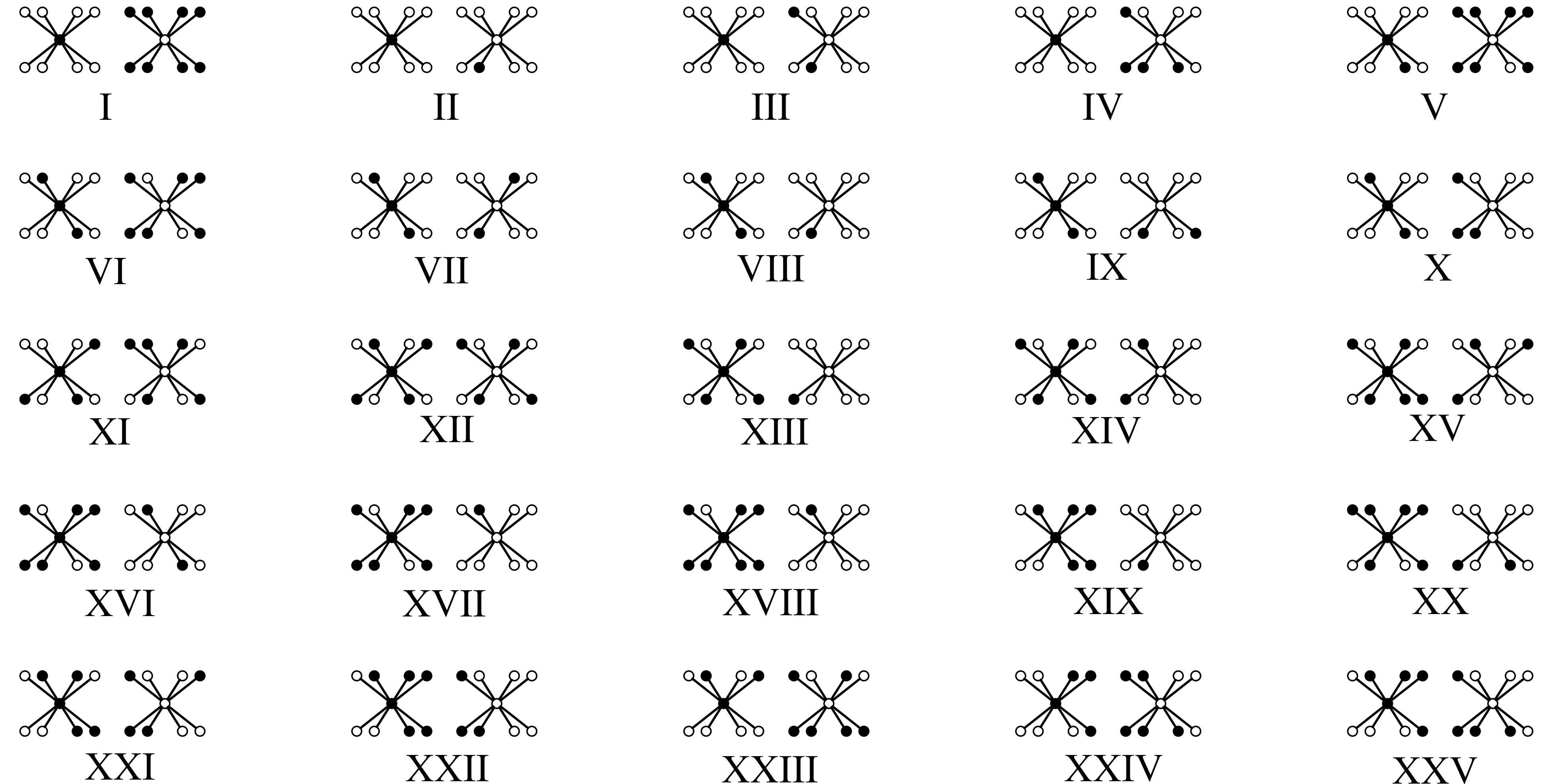}
\caption{Types of node neighbourhoods in the two-orbit partitions of $\Lambda^4$.}
\label{f:pcu4d}
\end{figure}
\end{center}

\newpage

{\tiny
\begin{table}[!h]
\begin{center}
\caption{Correspondences between two-orbit partitions of $\Lambda^4$ and $\Lambda^3$.}
\label{t:table6}
\begin{tabular}{|l|l|c|}
\hline
4D $Aut(\Pi)$ \{Stabilizers\} & 3D $Aut(\Pi)$ \{Stabilizers\} & Table~\ref{t:table4} \\ \hline
25/11/4/1 $\{C_2 \times C_2 \times S_4; C_2 \times C_2 \times S_4\}$ & $Fm\bar3m \, \{C_2 \times S_4; C_2 \times S_4\}$ & 1 \\
\hline
25/11/1/1 $\{C_2 \times C_2 \times S_4; C_2 \times C_2 \times S_4\}$ & $P4/mmm \, \{C_2 \times D_4; C_2 \times D_4\}$ & 3 \\
\hline
25/11/2/1 $\{C_2 \times C_2 \times S_4; C_2 \times C_2 \times D_4\}$ & $Im\bar3m \, \{C_2 \times S_4; C_2 \times D_4\}$ & 2 \\
\hline
25/11/1/1 $\{C_2 \times C_2 \times S_4; C_2 \times S_4\}$ & $P4/mmm \, \{C_2 \times D_4; D_4\}$ & 5 \\
\hline
19/6/1/1 $\{D_4 \times D_4; D_4 \times D_4\}$ & $P4/mmm \, \{C_2 \times D_4; C_2 \times D_4\}$ & 4 \\
\hline
25/11/1/1 $\{C_2 \times S_4; C_2 \times S_4\}$ & $P4/mmm \, \{D_4; D_4\}$ & 6 \\
\hline
13/10/2/1 $\{C_2 \times C_2 \times D_4; C_2 \times D_4\}$ & $Cmmm \, \{(C_2)^3; C_2 \times C_2\}$ & 11 \\
\hline
19/3/1/1 $\{C_4 \times D_4; D_4\}$ & $P4/m \, \{C_4 \times C_2; C_2\}$ & 18 \\
\hline
14/10/3/1 $\{C_2 \times D_6; D_6\}$ & $P\bar3m1 \, \{D_6; S_3\}$ & 8 \\
\hline
19/6/1/1 $\{C_2 \times D_4; C_2 \times D_4\}$ & $P4/mmm \, \{C_2 \times C_2; C_2 \times C_2\}$ & 12 \\
\hline
13/10/3/1 $\{C_2 \times D_4; C_2 \times D_4\}$ & $I4/mmm \, \{D_4; D_4\}$ & 7 \\
\hline
13/10/2/1 $\{C_2 \times D_4; C_2 \times D_4\}$ & $Cmmm \, \{C_2 \times C_2; C_2 \times C_2\}$ & 14 \\
\hline
13/10/1/22 $\{C_2 \times D_4; C_2 \times D_4\}$ & $Pmma \, \{C_2 \times C_2; C_2 \times C_2\}$ & 15 \\
\hline
14/4/1/1 $\{C_6 \times C_2; C_2\}$ & $R\bar3 \, \{C_6; C_1\}$ & 23 \\
\hline
25/11/4/1 $\{D_6; D_6\}$ & $Fm\bar3m \, \{S_3; S_3\}$ & 9 \\
\hline
25/11/2/1 $\{D_6; C_2 \times C_2\}$ & $Im\bar3m \, \{S_3; C_2\}$ & 19 \\
\hline
25/11/2/3 $\{D_6; C_2 \times C_2\}$ & $Ia\bar3d \, \{S_3; C_2\}$ & 20 \\
\hline
14/10/1/1 $\{D_6; D_6\}$ & $R\bar3m \, \{S_3; S_3\}$ & 10 \\
\hline
13/10/3/1 $\{(C_2)^3; (C_2)^3\}$ & $I4/mmm \, \{C_2 \times C_2; C_2 \times C_2\}$ & 13 \\
\hline
13/10/3/33 $\{(C_2)^3; (C_2)^3\}$ & $I4_1/amd \, \{C_2 \times C_2; C_2 \times C_2\}$ & 17 \\
\hline
6/3/2/3 $\{(C_2)^3; (C_2)^3\}$ & $Cmcm \, \{C_2 \times C_2; C_2 \times C_2\}$ & 16 \\
\hline
13/10/3/33 $\{C_2 \times C_2; C_2 \times C_2\}$ & $I4_1/amd \, \{C_2; C_2\}$ & 22 \\
\hline
6/3/2/3 $\{C_2 \times C_2; C_2 \times C_2\}$ & $Cmcm \, \{C_2; C_2\}$ & 21 \\
\hline
13/10/3/33 $\{C_2; C_2\}$ & $I4_1/amd \, \{C_1; C_1\}$ & 24 \\
\hline
13/10/3/37 $\{C_2; C_2\}$ & $I4_1/acd \, \{C_1; C_1\}$ & 25 \\
\hline
\end{tabular}
\end{center}
\end{table}
}

\newpage

The pioneering work of Heesch on two-orbit colourings remained absolutely unnoticed for almost one hundred years. In his authoritative biography \cite{Big88} the paper \cite{Hee33} is merely included in the list of publications but not commented upon in any way. By calling attention to this paper, we would like to acknowledge the contributions of H. Heesch to mathematical crystallography. Next year there will be the 120th anniversary of his birth, and this year (2025) marks the 30th anniversary since his death in 1995.

Heesch observed that two-orbit partitions of cubic lattices where unlike nodes have complementary neighbourhoods (\emph{cf.} Figs.~\ref{f:sql2}, \ref{f:pcu2}, \ref{f:pcu4d}) show global \emph{Verteilungssymmetrie} (permutation symmetry) in the sense that it does not matter which orbit is coloured black and which is white as in the case of the chessboard pattern. The numbers of such partitions are 6, 17, 62 for $d=2, 3, 4$, resp. This property is connected with the fact that two-orbit partitions represent a class of partitions with a high degree of symmetry. More specifically, every two-orbit partition into orbits of equal cardinality (hence the same number of unlike neighbours) in Tables~\ref{t:table4}, \ref{t:table5} arises from an index-2 subgroup $H$ of some node-transitive group $G$ of $\Lambda^d$. Indeed, $H$ is either again transitive or acts on $\Lambda^d$ with two orbits of equal cardinality which are mapped onto each other by elements from $G \setminus H$. It is possible to construct partitions of $\Lambda^d$ of a more general kind (for example, from a broader class of \emph{equitable partitions} \cite{Gods93}) that are aperiodic but where nodes have nevertheless equal numbers of unlike neighbours \cite{Puz08}. We note here that equitable partitions for the square lattice have been constructed for up to nine colours \cite{Krot08, Krot23}, whereas in higher dimensions ($d \geq 3$) no results are available.

In a recent publication chromatic numbers of four-dimensional lattices were considered \cite{Frank25}. Therein, chromatic numbers were determined for the Cayley graphs of $\Z^4$ corresponding to generating sets that consist of facet vectors of Voronoi parallelohedra. The chromatic number is defined as the minimum number of colours which are necessary to colour vertices of a graph in a way that no vertices of the same colour are adjacent. Our approach can be used to compute orbit colourings with the minimum number of colours. Cubic lattices considered in our work represent however the most transparent case since they give rise to bipartite graphs (their chromatic number is 2).

Orbit partitions (colourings) are relevant for coding theory, for certain problems in statistical physics \cite{Bax70}, for modelling distribution of atoms in solid solutions \emph{etc.} The theory of orbit colourings in connection with crystallographic structures (tilings, lattices, \emph{etc.}) was dealt with in a series of publications from the group of Prof. De las Pe\~nas and collaborators \cite{Jun23, Felix11}. The paper \cite{Felix11} contains a detailed exposition of orbit colourings for lattices together with a few examples for $\Lambda^2$ and $\Lambda^3$. An important result of \cite{Felix11} concerns the computation of the automorphism group for the orbit partition~$\Pi$ of $\Lambda^d$. In our computations $Aut(\Pi)$ was determined from ordering by inclusion the subgroups of $Aut(\Lambda^d)$ which act on $\Lambda^d$ with a given number of orbits. If, however, a comprehensive subgroup search is not performed, in \cite{Felix11} an elegant recipe is proposed for computing $Aut(\Pi)$. Suppose a group $H<Aut(\Lambda^d)$ acts on $\Lambda^d$, and consider the partition~$\Pi$ of $\Lambda^d$ into the orbits of~$H$. The $H$-orbits originate from the cosets of $T(H)$ in $\Lambda^d$. $Aut(\Pi)$ is found from the stabilizers of the $H$-orbits in the normalizer of $T(H)$ in $Aut(\Lambda^d)$. The computation of $Aut(\Pi)$ involves in general two steps that is not mentioned in \cite{Felix11}. Although $T(N_{Aut(\Lambda^d)}(T(H))) = \Lambda^d$ always holds, nevertheless the point group of $N_{Aut(\Lambda^d)}(T(H))$ might be too low to find all automorphisms. The intersection of stabilizers of the $H$-orbits in $N_{Aut(\Lambda^d)}(T(H))$ therefore gives rise to the intermediate group $S$. From the normalizer $N_{Aut(\Lambda^d)}(T(S))$ the full group $Aut(\Pi)$ is readily determined.

\medskip

\noindent \emph{Example:} The crystal structure of calcite $\rm{CaCO_3}$ (sp.gr. $R\bar3c$) is usually described in relation to the rocksalt type (NaCl) if the positions of Ca-atoms are identified with those of Na, and the barycenters of carbonate-anions -- with those of Cl. Moreover, cations and anions occupy Wyckoff positions with non-isomorphic stabilizers: Ca-atoms are situated at the sites with the $\bar3$ symmetry, anions -- at the sites with the $D_3$ symmetry. We have:

\begin{center}
\begin{tabular}{l}
$Aut(\Lambda^3) = \langle a, b, t_x \rangle \cong Pm\bar3m$, \\
$H = \langle a, m_{xxz}{t_y}^2, {t_x}^2t_yt_z \rangle \cong R\bar3c$, \\
$N_{Aut(\Lambda^3)}(T(H)) = \langle a, m_{xxz}, t_x \rangle \cong R\bar3m$, \\
\end{tabular}
\end{center}

\noindent where $a = \bar3_{xxx}: -z, -x, -y; b = m_{x\bar{x}z}: -y, -x, z; m_{xxz} = (aba)^2b$, and $t_x, t_y, t_z$ are unit translations in the respective directions. The $H$-orbits mod $T(H)$ are

\centerline{$\{ (0, 0, 0), (0, 0, 2) \}, \{ (0, 0, 1), (0, 0, 3) \}$.}

\noindent It is easy to see that $S = \langle a, m_{xxz}, t_yt_z \rangle$. The $S$-orbits mod $T(S)$ are $\{ (0, 0, 0) \}, \{ (0, 0, 1) \}$. Finally, we obtain

\centerline{ $N_{Aut(\Lambda^3)}(T(S)) = Aut(\Lambda^3); Aut(\Pi) = \langle a, b, t_yt_z \rangle \cong Fm\bar3m$.}

\medskip

\noindent Interestingly, there is only one more alternative to get the NaCl pattern with stabilizers $\bar3$ and $D_3$: it is realised in the subgroup $\langle a, a^3bt_x{t_y}^{-1}, aba^2t_x{t_z}^{-1} \rangle \cong Fd\bar3c$.
\bigskip

\noindent {\bf Acknowledgements}. The author thanks Prof. Denis Krotov for helpful correspondence on the partitions of the square lattice and for confirming the number of the 3-orbit partitions. Part of this work was done in 2022--2023 when the author served as a substitute professor (Vertretungsprofessor) at the Ludwig-Maximilians-Universit\"at M\"unchen, Sektion Kristallographie. Funding within the Hightech Agenda Bayern is gratefully acknowledged.

\bigskip

\noindent {\bf Supporting information}: Generating sets for node-transitive groups of cubic lattices (Table~\ref{t:table3}) and the data on partitions from Tables~\ref{t:table3}--\ref{t:table5} are available from the author upon request.

\end{document}